\documentclass[11pt]{article}
\usepackage[english]{babel}
\usepackage{graphics}
\usepackage{amsfonts}
\usepackage{amsmath}
\usepackage{amssymb}
\usepackage{mathrsfs}
\usepackage{enumerate}

\usepackage{hyperref}
\usepackage{epsfig}
\usepackage{amscd}
\usepackage{amsmath}

\newtheorem{lemma}{Lemma}[section]

\newtheorem{proposition}{Proposition}[section]

\newtheorem{theorem}{Theorem}[section]

\newtheorem{remark}{Remark}[section]
\newtheorem{definition}{Definition}[section]

\def\scatola{\lower5pt\hbox{\vbox{\hrule\hbox{\vrule\kern2pt\vbox%
{\kern5pt\hbox{\mathsurround=0pt }\kern2pt}\kern4pt\vrule}\hrule}}\
} 
\def\l{\lambda}

\def\H{{I\!\!H}}

\def\R{{\mathbb{R}}}

\def\l{\lambda}

\def\0{\overline 0}

\def\Hcal{\mathcal{H}}

\def\Hca{{\cal H}}

\font\tenmsb=msbm10 \font\sevenmsb=msbm7 \font\fivemsb=msbm5
\newfam\msbfam
\textfont\msbfam=\tenmsb \scriptfont\msbfam=\sevenmsb
\scriptscriptfont\msbfam=\fivemsb
\font\teneufm=eufm10 \font\seveneufm=eufm7 \font\fiveeufm=eufm5
\newfam\eufmfam
\textfont\eufmfam=\teneufm \scriptfont\eufmfam=\seveneufm
\scriptscriptfont\eufmfam=\fiveeufm
\def\frak#1{{\fam\eufmfam\relax#1}}

\newcommand*{\finedim}{{\hfill{\lower5pt\hbox{\vbox{\hrule\hbox{\vrule\kern2pt\vbox%
{\kern5pt\hbox{\mathsurround=0pt
}\kern2pt}\kern4pt\vrule}\hrule}}\bigskip }}}
\newcommand*{\de}{\,\mathrm{d}}
\newcommand*{\xp}{{\dot x}}

\begin{document}

\title{Differential games and Hamilton--Jacobi equations in the Heisenberg group}
\author{ A. Calogero\thanks{
Dipartimento di Matematica e Applicazioni, Universit\`a degli Studi
di Milano--Bicocca, Via Cozzi 53, 20125 Milano, Italy ({\tt
andrea.calogero@unimib.it}) }}

\date{}
\maketitle

\begin{abstract}
\noindent
The purpose of this work is twofold.
First we study the solutions of a Hamilton--Jacobi equation of the form $u_t(t,x)+\Hcal(t,x,\nabla_H u(t,x))=0$, where $\nabla_H u$ represents the horizontal gradient of a function $u$ defined on the Heisenberg group $\H$. Motivated by \cite{LiMaZh2016}, we prove a  Lipschitz continuity preserving property  for  $u$ with respect to the  Kor\'anyi homogeneous distances $d_G$ in $\H$.    Secondly, we are keenly interested in introducing the game theory in $\H$, taking into account its Sub--Riemannian structure: inspired by \cite{EvSo1984} and \cite{BaCaPi2014}, we prove $d_G$-Lipschitz regularity results for the lower and the upper value functions of a zero game with horizontal curves as its trajectories, and we study the Hamilton--Jacobi--Isaacs equations associated to such zero game. As a consequence, we also provide a representation of the viscosity solution of the initial Hamilton--Jacobi equation.

\end{abstract}


\medskip



\section{Introduction}

The first aim of this paper is to study the properties of the viscosity solutions of the  Hamilton--Jacobi equation in the Heisenberg group $\H$
\begin{equation}\label{Ham Jac eq intro}
\left\{
  \begin{array}{ll}
    \displaystyle \frac{\partial u}{\partial t}(t,x)+\Hcal(t,x,\nabla_{H}u(t,x))=0 & \hbox{in $(0,T)\times\H$} \\
    u(0,x)=g(x) & \hbox{in $\H$,}
  \end{array}
\right.
\end{equation}
where  $\nabla_H u$ is the horizontal gradient of the function $u$, and  $\Hcal$ and $g$ are bounded  functions satisfying suitable assumptions (see \emph{3.} and \emph{4.}); in particular, they are
Lipschitz continuous in $x$ w.r.t. the left--invariant Kor\'anyi distance $d_G$.

This study is motivated by the work of Liu, Manfredi and Zhou \cite{LiMaZh2016}; however, our approach is different. It is well known that the study of the  Hamilton--Jacobi equations is strictly related to the game theory; the pioneer of this approach was Isaacs \cite{Is1965}. Several authors like Evans, Souganidis, Bardi, Lions   have connected the Isaacs theory with the notion of viscosity solution (see \cite{EvSo1984}, \cite{BaEv1984} \cite{LiSo1985}).
Along this line of investigation, and since we are  interested in embedded the game theory in the Sub--Riemannian structure of the Heisenberg group $\H$, we study the viscosity solution of \eqref{Ham Jac eq intro} and its properties by introducing the following  zero game:
\begin{equation}\label{zero sum intro}
\left\{
  \begin{array}{ll}
    \displaystyle{\texttt{\rm Player I:}\ \max_{y\in {\mathcal Y}(0)} J(y,z)\qquad\qquad \texttt{\rm Player II:}\ \min_{z\in {\mathcal Z}(0)} J(y,z)} & \\
    \displaystyle{J(y,z)
=\int_0^T F(t,x,y,z)\de t+g(x(T))} & \\
x:[0,T]\to \H\ \ {\texttt{\rm  horizontal curve with }}\ z=(\dot{x}_1,\dot{x}_1)\ {\texttt{\rm  a.e. }}&\\
x(0)=x_0,&\\
  \end{array}
\right.
\end{equation}
where $T>0$ and $x_0\in\H$ are fixed, the function $F$ satisfies assumption \emph{2.} and the function $g$ is as in \eqref{Ham Jac eq intro}; in particular,
such two functions  are
Lipschitz continuous in $x$ w.r.t. $d_G$.
The set of controls for the two players are  defined by
${\mathcal Y}(0)=\{ y:[0,T]\to Y\subset\R^2,\ \texttt{\rm measurable}\}$ and
 ${\mathcal Z}(0)=\{
z:[0,T]\to Z\subset\R^2,\ \texttt{\rm measurable}\}.$

This idea to study Hamilton--Jacobi equations using optimal control theory
 was successfully followed in \cite{BaCaPi2014}, where the generic  dynamics is replaced by the condition for $x$ to be a horizontal curve on $\H$ (see subsection \ref{curve orizzonatali} for details).
 However, to our knowledge of the literature, this is the first work that connects game theory and Heisenberg group; this connection is the second aim of this paper.

In order to construct such connection between Heisenberg group, game theory and Hamilton--Jacobi equation, we introduce the following, and  classical, assumptions for the problem
\eqref{Ham Jac eq intro}
and for the game
 \eqref{zero sum intro}:

\begin{itemize}
\item[\emph{1.}] the control sets $Y$ and $Z$ are compacts; more precisely, $Y=\{y\in\R^2: \ \|y\|\le R_Y\},$  $Z=\{z\in\R^2: \ \|z\|\le R_Z\}$,  for some fixed and positive $R_Y$ and $R_Z$;
  \item[\emph{2.}] the function $F:[0,T]\times\H\times Y\times Z\to\R$ is uniformly continuous, with
$$
\left\{
  \begin{array}{ll}
    |F(t,x,y,z)|\le C_1 &  \\
 |F(t,x,y,z)- F(t,x',y,z)|\le C_1' d_G(x,x') &  \end{array}
\right.
$$
for some constants $C_1,\ C_1'$ and for every $ t\in [0,T]$,  $x,x'\in\H$, $y\in Y$ and $z\in Z$;
  \item[\emph{3.}] the function $g:\H\to\R$ satisfies
$$\left\{
  \begin{array}{ll}
    |g(x)|\le C_2 &  \\
  |g(x)- g(x')|\le C_2' d_G(x,x')&
  \end{array}
\right.
$$
for some constants $ C_2,\ C_2'$ and for every   $x,x'\in\H$;

\item[\emph{4.}] the function  $\Hcal:[0,T]\times\H\times Y\to\R$ is uniformly continuous, with
$$
\left\{
  \begin{array}{ll}
    |\Hcal(t,x,y)|\le D_1 &  \\
 |\Hcal(t,x,y)- \Hcal(t,x',y)|\le D_1' d_G(x,x') &\\
|\Hcal (t,x,y)- \Hcal(t,x,y')|\le K \|y-y'\|&  \end{array}
\right.
$$
for some constants $D_1, D_1', K,$ and for every $ t\in [0,T]$,  $x,x'\in\H$, $y, y'\in \R^2.$
\end{itemize}

\noindent
We emphasize that in assumptions \emph{2.}, \emph{3.}  and  \emph{4.}  a $d_G$-Lipschitz condition is considered, where the mentioned Kor\'anyi--gauge metric $d_G$ is a natural metric in $\H$ that turns out to be equivalent to the  Carnot--Carath\'eodory metric $d_{CC}$. If we replace these $d_G$-Lipschitz properties  with the Euclidean Lipschitz requirements, we know that the lower $V^-$ and the upper $V^+$ value functions for the game \eqref{zero sum intro} (see Definition \ref{def upper and lower value functions})
are Euclidean Lipschitz (see Remark \ref{euclidean assumption}). Our result is indeed more precise and in the spirit of proving Lipschitz preserving properties as in \cite{LiMaZh2016}:

\begin{theorem}[$d_G$-Lipschitz continuity preserving properties for $V^-$]\label{Teo visco}

Let us consider the zero game (\ref{zero sum intro}) with the assumptions 1., 2.  and 3..
 Then its lower value function $V^-$ is bounded and uniformly
Lipschitz continuous w.r.t. the metric $d_G$, i.e. there exists a constant $C'$ such that
$$|V^-(t,x)-V^-(t',x')|\le  C'\left(|t-t'|+d_G(x,x')\right),$$
for every $t,t'\in [0,T]$ and $x,x'\in\H$.  A similar result holds for $V^+$.
 \end{theorem}

The second important result of the paper is the following:

\begin{theorem}[$V^-$ as viscosity solution]\label{Teo vis2 intro}
Let us consider the zero game (\ref{zero sum intro}) with the assumptions
1., 2. and 3.. Then, the lower value function $V^-$ is a viscosity solution
of the lower Hamilton--Jacobi--Isaacs equation
 \begin{equation}\label{system PD-}
 \left\{
   \begin{array}{ll}
   {\displaystyle \frac{\partial u}{\partial t}(t,x)+H^-(t,x,\nabla_{H}u(t,x))=0} & \hbox{for $(t,x)\in(0,T)\times\H$} \\
  u(T,x)=g(x)&  \hbox{for
$x\in\H$,}
   \end{array}
 \right.
 \end{equation}
 where  $H^-$ is the lower Hamiltonian for the game \eqref{zero sum intro}.
  \end{theorem}
 The  lower Hamiltonian $H^-$  is defined as a maxmin-function (see Definition \ref{def upper lower ham}).
Clearly, a similar result holds for $V^+$.
The proofs of Theorem \ref{Teo visco} and Theorem \ref{Teo vis2 intro} require a fine use of the horizontal curves and their properties in $\H$.

A precise estimate of the $d_G$-Lipschitz constant for $V^-$ in Theorem \ref{Teo visco} allows us to provide a representation of the viscosity solution $u$ for the initial problem \eqref{Ham Jac eq intro} as a value function (see Theorem \ref{ciao}): as a consequence of the previous results, the mentioned $d_G$-Lipschitz assumptions in \emph{3.} and \emph{4.} for the functions $\Hcal$ and $g$ involved in \eqref{Ham Jac eq intro}
is inherited by $u$, that turns out to be $d_G$-Lipschitz.

The paper is organized as follows: in Section 2 we introduce and connect the fundamental notions in the Heisenberg group and in the game theory; moreover, we prove some fine properties of the horizonal curves in $\H$. Section 3 and Section 4 are essentially devoted to the proofs of Theorem \ref{Teo visco}
 and Theorem \ref{Teo vis2 intro}, respectively. In Section 5 we study  problem \eqref{Ham Jac eq intro} and raise an open question related to the Hopf--Lax formula.


\section{Preliminaries.}

\subsection{A short introduction on the Heisenberg group $\H$}

The Heisenberg group $\H$ is $\R^3$ endowed with a non--commutative law $\circ$: it is the Lie group whose
 Lie algebra $\frak{h}$ admits a stratification of step 2; in particular
 $\frak{h}=\mathbb{R}^3=V_1\oplus V_2,$
with
\begin{equation}\label{algebra Heisenberg}
\begin{array}{ll}
V_1={\rm span}\left\{X_1,X_2\right\}\quad&\texttt{\rm with}\ \
X_1=\partial_{x_1}-\frac{x_2}{2}\partial_{x_3}\quad \texttt{\rm
and}\ \
X_2=\partial_{x_2}+\frac{x_1}{2}\partial_{x_3},\\
  V_2={\rm
span}\left\{X_3\right\}&\texttt{\rm with}\ \ X_3=\partial_{x_3}.
\end{array}
\end{equation}
The bracket $[\cdot ,\cdot ]:\frak{h} \times\frak{h} \to \frak{h}$
is defined as $[X_1,X_2]=X_3,$ and it vanishes for all the other
basis vectors.
For a sufficiently regular function $f:\H\to\R$, we define the horizontal gradient $\nabla_H f$ by
$$\nabla_H f(x)=\left( X_1 f(x),\ X_2 f(x)\right)\qquad \texttt{\rm with}\ x\in\H ;$$
 for our purpose it is convenient to think  $\nabla_H f(x)$ as a vector in $\R^2$.
We say that such $f$ is in $\Gamma^1(\H)$ if its horizontal derivatives $ X_1 f$ and $X_2 f$ are continuous functions.

 The group law is defined by the relation
$$
(x_1,x_2,x_3)\circ(x_1',x_2',x_3')=
(x_1+x_1',x_2+x_2',x_3+x_3'+(x_1x_2'-x_1'x_2)/2).$$
Consequently, the null element is $e=(0,0,0)$ and $(x_1,x_2,x_3)^{-1}=(-x_1,-x_2,-x_3).$ The dilation is a family of automorphisms given
by $\delta_\lambda(x_1,x_2,x_3)=(\lambda x_1,\lambda x_2,\lambda^2
x_3),$ and hence the homogeneous dimension is 4.

We denote by
$\|\cdot\|$ the Euclidean norm in $\R^n$ and by
$d_E$ the Euclidean distance. In contrast with Analysis
in Euclidean spaces, where the Euclidean distance is the most natural
choice, in the Heisenberg group several distances have been introduced for
different purposes (for example, see $d_{CC}$ below).
The Kor\'anyi gauge  $\|\cdot\|_G$, defined  by
\begin{equation}\label{gauge}
 \| (x_1,x_2,x_3)\|_G=\left((x_1^2+x_2^2)^2+x_3^2\right)^{1/4}\qquad \forall x=(x_1,x_2,x_3)\in\H,
\end{equation}
allows us to introduce the left invariant metric
 $d_G$  via  $d_G(x,x')=\|(x')^{-1}\circ x\|_G.$
This  Kor\'anyi distance $d_G$ is  homogeneous, namely,
it is continuous, left invariant and behaves well respect to the dilations $\delta_\lambda$.
For our purpose it is important to mention that
for every compact set $\Omega\subset\H$ there exists
a constant $C=C(\Omega)$ such that (see for example \cite{BoLa2012})
\begin{equation}\label{dis euclidea e gauge}
\frac{1}{C}\|x\|\le  \|x\|_G\le C \|x\|^{1/2}\qquad \forall x\in\Omega.
\end{equation}

Given a metric $d$  in $\H$, we say that a function $f:\H\to \R$ is Lipschitz w.r.t. the metric $d$ (or shortly is $d$-Lipschitz) it there exists a constant $C$ such that
$$
|f(x)-f(x')|\le C d(x,x')\qquad\forall x,x'\in \H.
$$
In order to emphasize the dependence on the distance $d$, we stress that the function $f(x)=\|x\|_G$ is $d_G$-Lipschitz, but it is not $d_E$-Lipschitz. On the other hand, by \eqref{dis euclidea e gauge}, it is clear that every  $d_E$-Lipschitz function is  $d_G$-Lipschitz.

\subsection{Horizontal curves in $\H$}\label{curve orizzonatali}

Let us start with this fundamental notion:
\begin{definition}[horizontal curve]
A horizontal curve
$x=(x_1,x_2,x_3):[a,b]\to\H$, with $[a,b]\subset\R$, is an absolutely continuous function a.e. tangent
to horizontal directions, i.e.
$$\dot{x}(s)=\dot{x}_1(s) X_1(x(s))+\dot{x}_2(s) X_2(x(s))\qquad \texttt{\rm a.e.}\ s\in [a,b].
$$
\end{definition}
Equivalently, $x$ is horizontal if
$$
\dot{x}_3=(x_1\dot{x}_2-x_2\dot{x}_1)/2\qquad \texttt{\rm a.e.\ in}\ [a,b],
$$
that is  $\dot{x}=f^\H(x,z)$, for some measurable function $z:\R\to\R^2$, where
$f^\H:\R^3\times\R^2\to\R^3$ is given by
\begin{equation}\label{dynamic heisenberg}
f^\H(x,z)=\left(\begin{array}{c}
z_1\\
z_2\\
(z_2x_1-z_1x_2)/2
\end{array}\right)\qquad \forall z\in\R^2, \ x\in\H.
\end{equation}
We say that $x$ is a horizontal curve on $[a,b]\subset\R$ with horizontal velocity $z$ and initial point $\xi\in\H$ if
$$
x(t)=\xi+\int_a^t f^\H(x(s),z(s))\de s\qquad \forall s\in[a,b].
$$
The apex \lq\lq$\H$\rq\rq\ in \eqref{dynamic heisenberg} for the function $f$ reminds that, hereinafter and, in particular, in the framework of our games,
we are deal with horizontal curves in the Heisenberg Sub--Riemannian geometry.

In the sequel we will often use the dynamics   $\dot{x}=-f^\H(x,z)$ for a horizontal curve: it is easy to see that, for a given  $z$, the curves $x$ and $\widetilde x$ defined by
\begin{equation}\label{curve hor meno}
\left\{
  \begin{array}{ll}
   \xp=-f^\H( x,z)\qquad {\texttt{\rm  a.e. in}}\ [a,b]&\\
x(a)=\xi&\\
  \end{array}
\right.
\qquad\qquad
\left\{
  \begin{array}{ll}
   \dot{\widetilde x} =f^\H( \widetilde x,- z)\qquad {\texttt{\rm  a.e. in}}\ [a,b]&\\
\widetilde x(a)=\xi&\\
  \end{array}
\right.
\end{equation}
 are the same object, i.e. the horizontal curve on $[a,b]$ with velocity $-z$ and initial point $\xi$.

Given $\xi$ and $\xi'$ in $\H$, we denote by $\Gamma(\xi,\xi')$ the set of all horizontal curves $x$ in $[0,1]$ with initial point $x(0)=\xi$ and final point $x(1)=\xi'$. Chow's theorem (see, for example, \cite{CaDaPaTy2007}) guarantees that $\Gamma(\xi,\xi')\not=\emptyset$ and allows us to introduce the
 Carnot--Carath\'eodory metric $d_{CC}$ by the formula
 $$
d_{CC}(\xi,\xi')=\inf_{x\in\Gamma(\xi,\xi')}\int_0^1 \|(\dot{x}_1(s),\dot{x}_2(s))\|\de s.
 $$
 It is well--known that the Carnot--Carath\'eodory metric $d_{CC}$ and the gauge metric $d_G$ are bi--Lipschitz equivalent.

The following property is known and can be easily proved:

\begin{remark}\label{derivata} Let $F:\H\to\R$ be a function in $\Gamma^1$, and let
 $x=(x_1,x_2,x_3):[0,T]\to\H$ be a  horizontal curve.
 Then, for a.e. $s\in (0,T)$,
$$
\frac{d F(x(s))}{d s}=(\dot x_1(s),\dot x_2(s))\cdot \nabla_H F(x(s)).
$$
\end{remark}

The next three propositions are crucial in order to prove our Lipschitz preserving property for the value functions.
For every fixed $\tau\in[0,T]$,  let us introduce the set of controls at
time $\tau$ for Player I as
 ${\mathcal Z}(\tau)=\{
z:[\tau,T]\to Z\subset\R^2,\ \texttt{\rm measurable}\}.$

\begin{proposition}\label{bounded traj} Let assumption \emph{1.} be satisfied.
Let $\tau$ be fixed in $[0,T]$, and let us consider a horizontal curve
$x=(x_1,x_2,x_3):[\tau,T]\to\H$ with horizontal velocity  $z=(z_1,z_2)\in {\mathcal Z}(\tau)$ and initial point $\xi\in\H$, i.e.
\begin{equation}\label{curve hor 0}
\left\{
  \begin{array}{ll}
      \dot{x}=f^\H  (x,z) \qquad \texttt{\sl a.e. in } [\tau,T] &\\
      x(\tau)=\xi &\\
  \end{array}
\right.
\end{equation}
Then
\begin{equation}\label{bounded traj1} d_G(\xi,x(t))\le  3 R_Z(t-\tau)\qquad\forall t\in[\tau,T]\end{equation}
\end{proposition}
\textbf{Proof} The assertion \eqref{bounded traj1} comes from
\begin{eqnarray*}
\left|\left(\xi^{-1}\circ x(t)\right)_i\right|&\le&|x_i(t)-\xi_i|=\int_{\tau}^t|z_i(s)|ds\le (t-\tau)R_Z,\qquad i=1,2\\
\left|\left(\xi^{-1}\circ x(t)\right)_3\right|&=&
\left|-\xi_3+x_3(t)+\frac{1}{2}\left(-\xi_1 x_2(t)+\xi_2 x_1(t)\right)\right|\\
&=&\frac{1}{2}
\left|\int_{\tau}^t\left[z_2(s)(x_1(s)-\xi_1)-z_1(s)(x_2(s)-\xi_2)\right]\de s\right|\\
&\le&
\frac{1}{2}\int_{\tau}^t\left|(x_1(s)-\xi_1,x_2(s)-\xi_2)\cdot(z_2(s),-z_1(s))\right|\de s\\
&\le &(t-\tau)^2R_Z^2.
\end{eqnarray*}
Since $\|\xi^{-1}\circ x(t)\|_G=d_G(\xi, x(t))$, by \eqref{gauge} we have the claim.\finedim

\begin{proposition}\label{hor curve new 1}
 Let us suppose that assumption 1. holds, and let $x$ be a horizontal curve as in Proposition \ref{bounded traj}.
Let $\widehat\xi$ be fixed in $\H$ and let us consider the horizontal curve $\widehat x$ on $[\tau,T]$ with horizontal velocity $z$, i.e.
\begin{equation}\label{curve hor}
\left\{
    \begin{array}{ll}
      \dot{\widehat x}=f^\H  (\widehat x,z)& \hbox{in $[\tau,T]$} \\
      \widehat x(\tau)=\widehat\xi &     \end{array}
  \right.
\end{equation}
Then there exists a constant $\widehat C $ that depends only on $R_Z$ and $T$ such that
\begin{equation}\label{curve hor estimate 1}
d_G(x(t),\widehat x (t))\le \widehat C d_G(\xi,\widehat\xi)\qquad\qquad \forall t\in [\tau,T].
\end{equation}
\end{proposition}
\textbf{Proof}
 Let us define $\phi:[\tau,T]\to [0,\infty)$ by $\phi(t)=d_G\left(x(t),\widehat x (t)\right)$. Then
\begin{eqnarray*}
\phi(t)&=&\left\|(\widehat x(t))^{-1}\circ x(t)\right\|_G\\
&=&\left\|\left(-\widehat x_1(t)+x_1(t),\ -\widehat x_2(t)+x_2(t),\ -\widehat x_3(t)+x_3(t)-\frac{1}{2}\left(\widehat x_1(t)x_2(t)-\widehat x_2(t)x_1(t)\right)\right) \right\|_G\\
&=&\Biggl\|\Biggl(-\widehat\xi_1-\int_\tau^t z_1(s)\de s+\xi_1+\int_\tau^t z_1(s)\de s,\\
&&\quad\ -\widehat\xi_2-\int_\tau^t z_2(s)\de s+\xi_2+\int_\tau^t z_2(s)\de s,\\
&&\quad\ -\widehat\xi_3-\frac{1}{2}\int_\tau^t  (\widehat x_1(s)z_2(s)-\widehat x_2(s)z_1(s))\de s
+\xi_3+\frac{1}{2}\int_\tau^t  ( x_1(s)z_2(s)- x_2(s)z_1(s))\de s+\\
&&\quad\qquad-\frac{1}{2}\Biggl[\left(\widehat\xi_1+\int_\tau^t z_1(s)\de s\right)\left(\xi_2+\int_\tau^t z_2(s)\de s\right)+\\
&&\qquad\qquad\qquad\  -
\left(\widehat\xi_2+\int_\tau^t z_2(s)\de s\right)\left(\xi_1+\int_\tau^t z_1(s)\de s\right)\Biggr]
   \Biggr)\Biggr\|_G\\
&=&\left\|\left(\xi_1-\widehat\xi_1,\
\xi_2-\widehat\xi_2,\
\xi_3 -\widehat\xi_3+\frac{1}{2}\left(\xi_1\widehat\xi_2-\xi_2\widehat\xi_1\right)+\int_\tau^t  ((\xi_1-\widehat \xi_1)z_2(s)-(\xi_2-\widehat \xi_2)z_1(s))\de s
   \right)\right\|_G
\end{eqnarray*}
By  \eqref{gauge} we have, a.e.,
\begin{eqnarray*}
\frac{d\phi(t)}{d t}
 &=&\frac{\left((\widehat x(t))^{-1}\circ x(t)\right)_3
 \left((\xi_1-\widehat \xi_1)z_2(t)-(\xi_2-\widehat \xi_2)z_1(t) \right)}{2\left(\phi(t)\right)^{3}}\\
 &=&
 \frac{\left|\left((\widehat x(t))^{-1}\circ x(t)\right)_3\right|}{2\left(\phi(t)\right)^{2}}
 \frac{\left|(\xi_1-\widehat \xi_1,\xi_2-\widehat \xi_2 )(z_2(t),-z_1(t) )\right|}{\phi(t)}\\
 &\le&\frac{1}{2}
\frac{\left\|\left(\xi_1-\widehat\xi_1,\ \xi_2-\widehat\xi_2\right) \right\|\left\|z(t)\right\|}{\left\|\left(\xi_1-\widehat\xi_1,\ \xi_2-\widehat\xi_2\right) \right\|}\\
&\le& \frac{R_Z}{2}
\end{eqnarray*}
Now, by the Gronwall inequality, we obtain
\begin{equation}\label{C hat 0}
\phi(t)\le \phi(\tau)\exp\left(\int_\tau^t  \frac{R_Z}{2} ds \right)\le \exp\left(T\frac{R_Z}{2}\right) d_G(\xi,\widehat \xi):=
\widehat C d_G(\xi,\widehat \xi).
\end{equation}
\vskip-0.4truecm\finedim

\noindent
Let us note that the curve $\widehat x$ in \eqref{curve hor} is exactly a left translation of the first curve $x$, i.e. $\widehat x(t)=\widehat\xi\circ\xi^{-1}\circ x(t).$

\begin{proposition}\label{hor curve new 2}
 Let us suppose that assumption 1. is satisfied and let $x$ be a horizontal curve as in Proposition \ref{bounded traj}.
Let    $\widetilde\xi$ in $\H$ and $0\le \tau\le \tau'\le T$ be fixed.
Let $\widetilde  x$ be the horizontal curve on $[\tau',T]$ with horizontal velocity $z\bigl|_{[\tau',T]}\in {\mathcal Z}(\tau')$ and initial point $\widetilde\xi$, i.e.
$$\left\{
    \begin{array}{ll}
      \dot{\widetilde x}=f^\H  (\widetilde x,z)& \hbox{in $[\tau',T]$} \\
      \widetilde x(\tau')=\widetilde\xi. &     \end{array}
  \right.
$$
Then there exists a constant $\widetilde C $ that depends only on $R_Z$ and $T$ such that
$$
d_G(x(t),\widetilde x (t))\le  \widetilde C\left(d_G(\widetilde\xi,\xi)+(\tau'-\tau)\right)\qquad\qquad \forall t\in [\tau',T].
$$
\end{proposition}
\textbf{Proof}
By Proposition \ref{bounded traj}, we have $d_G(\xi, x(\tau'))\le 3R_Z(\tau'-\tau).$
Now, it is easy to see that $\widetilde x$ is a left translation of the curve $x$ restricted to $[\tau',T]$: more precisely,
$\widetilde x(t)=\widetilde\xi\circ(x(\tau'))^{-1}\circ x(t)$.
By (\ref{curve hor estimate 1}) we have, for every $t\in[\tau',T]$,
\begin{eqnarray}
d_G\left(x(t),\widetilde x(t)\right)
 &\le&\widehat C  d_G\left(x(\tau'),\widetilde x(\tau')\right)\nonumber\\
  &\le& \widehat C \left(d_G(x(\tau'),\xi)+d_G(\xi, \widetilde\xi)\right)\nonumber\\
 &\le&  \widehat C (1+3R_Z)\left(d_G(\widetilde\xi,\xi)+(\tau'-\tau)\right)\nonumber\\
   &:=& \widetilde C\left(d_G(\widetilde\xi,\xi)+(\tau'-\tau)\right) \label{C tilde 0}
\end{eqnarray}
\vskip-0.4truecm\finedim

\subsection{Differential games in $\H$}

The game which we are interested in is the following (see \eqref{zero sum intro})
\begin{equation}\label{zero sum}
\left\{
  \begin{array}{ll}
    \displaystyle{\texttt{\rm Player I:}\ \max_{y\in {\mathcal Y}(0)} J(y,z)\qquad \qquad\texttt{\rm Player II:}\ \min_{z\in {\mathcal Z}(0)} J(y,z)} & \\
    \displaystyle{J(y,z)
=\int_0^T F(t,x,y,z)\de t+g(x(T))} & \\
\xp=-f^\H( x,z)\qquad {\texttt{\rm  a.e. in}}\ [0,T]&\\
x(0)=x_0&\\
  \end{array}
\right.
\end{equation}
under the
 assumptions\emph{ 1.}, \emph{2.} and \emph{3.},  where $T>0$ and $x_0\in\H$ are fixed.

Three comments on such zero game are required. Following the idea in \cite{BaCaPi2014}, the dynamics we consider is horizontal: essentially, we  consider a game where, for every strategy of the two players, the associated trajectory is a horizontal curve on $\H$.
Secondly, in the dynamics appears a minus (see \eqref{curve hor meno}) whose only reason is to keep consistency with the classical case. Finally, the reader who is not expert in game theory would be surprise by the fact that  the dynamics does not involve the control of the first Player. We want to reassure these readers, because this is a classical situation, asymmetric for the two players, and it turn out to be very useful to obtain representations for the solutions of  Hamilton--Jacobi equations.

Starting from this game, let us introduce the classical notions of controls and strategies for the two players (see for example \cite{BaCa1997}, \cite{BaOl1998}).

  For every fixed $\tau\in[0,T]$,  let us introduce the set of controls at
time $\tau$ for Player I as
${\mathcal Y}(\tau)=\{ y:[\tau,T]\to Y\subset\R^2,\ \texttt{\rm measurable}\}$.
In a similar way for the second Player, we define ${\mathcal Z}(\tau)=\{
z:[\tau,T]\to Z\subset\R^2,\ \texttt{\rm measurable}\}.$

\noindent
 We say that a map $\alpha:
{\mathcal Z}(\tau)\to{\mathcal Y}(\tau)$ is a
nonanticipative strategy for  Player I at time $\tau$
if, for any time $t\in[\tau,T]$ and any controls $z, z'\in
{\mathcal Z}(\tau)$ such that $z=z'$ a.e. in
$[\tau,t]$, then we have $\alpha[z] = \alpha[z']$ a.e. in $[\tau,t]$.
We denote by ${\mathcal S}_\alpha(\tau)$ the set of such
nonanticipative strategies at time $\tau$ for  Player
I.
 In a symmetric way, we denote by ${\mathcal
S}_\beta(\tau)$ the set of  nonanticipative strategies for Player II, which
are the nonanticipative maps
$\beta :
{\mathcal Y}(\tau)\to{\mathcal Z}(\tau)$.

\begin{definition}[upper and lower value functions]\label{def upper and lower value functions}
 Let us consider the zero game \eqref{zero sum}.   The lower  value function $V^-:[0,T]\times\H\to\R$ is defined by
\begin{equation}\label{def V-}
V^-(\tau,\xi)=\inf_{\beta\in {\mathcal S}_\beta(\tau)}\sup_{y\in
{\mathcal Y}(\tau)}\left\{ \int_\tau^T F(t,x,y,\beta[y])\de
t+g(x(T)) \right\} ,
\end{equation}
where $x$ is the horizontal curve on $[\tau,T]$ with horizontal velocity $-\beta[y]$ and initial point $\xi$.
   The upper  value function $V^+:[0,T]\times\H\to\R$ is defined by
\begin{equation}\label{def V+}
V^+(\tau,\xi)=\sup_{\alpha\in {\mathcal S}_\alpha(\tau)}\inf_{z\in
{\mathcal Z}(\tau)}\left\{ \int_\tau^T F(t,x,\alpha[z],z)\de
t+g(x(T)) \right\} ,
\end{equation}
where $x$ is the horizontal curve on $[\tau,T]$ with horizontal velocity $-z$ and initial point $\xi$.
\end{definition}
It is well know that in general $V^-\le V^+$ and such two functions are different (see \cite{BaCa1997} for an example of game where the previous inequality is strict). We say that the game (\ref{zero sum}) admits value function $V$ if
$$
V(\tau,\xi)=V^+(\tau,\xi)=V^-(\tau,\xi),\qquad\forall (\tau,\xi)\in [0,T]\times\H.
$$

\noindent
The following Dynamic Programming optimality condition is a classical result proved in \cite{ElKa1974}:
\begin{theorem}\label{Teo vis0}
Let us consider the problem (\ref{zero sum}).
 Then
$$
V^-(\tau,\xi)=\inf_{\beta\in {\mathcal S}_\beta(\tau)}\sup_{y\in
{\mathcal Y}(\tau)}\left\{ \int_\tau^{\tau+\sigma} F(t,x,y,\beta[y])\de
t+V^-(\tau+\sigma,x(\tau+\sigma)) \right\}
$$
for every $\tau,\ \tau+\sigma\in [0,T]$ and $\xi\in\H$. A similar result holds for $V^+$.
 \end{theorem}

\section{Lipschitz continuity preserving properties}

This section in devoted to the proof of Theorem \ref{Teo visco}.
First of all, let us remark that, if we consider the dynamics $\dot x=-f^\H(x,z)$ in game \eqref{zero sum}, under assumption \emph{1.}, it is easy to see that $f^\H$ is uniformly continuous with
$$
\|f^\H(x,z)-f^\H(x',z)\|=\frac{1}{2}\left|(z_2,-z_1)\cdot(x_1-x_1',x_2-x_2')\right|\le\frac{1}{2}R_Z\|x-x'\|,
$$
for all $x=(x_1,x_2,x_3),\ x'=(x'_1,x'_2,x'_3)$ in $\H$ and $z=(z_1,z_2)$ in $Z$. 
Now, let us replace in assumptions \emph{2.}  and \emph{3.}  the gauge distance $d_G$ with the Euclidean distance $d_E$, i.e. let us assume for a moment that
\begin{itemize}
  \item[\emph{2'.}] \qquad $|F(t,x,y,z)|\le C_1$,\qquad $|F(t,x,y,z)- F(t,x',y,z)|\le C_1' \|x-x'\|$

  \item[\emph{3'.}] \qquad $|g(x)|\le C_2$,\qquad \qquad $|g(x)- g(x')|\le C_2' \|x-x'\|$
\end{itemize}
for some constants $C_1,\ C_1',\ C_2,\ C_2'$ and for every $ t\in [0,T]$,  $x,x'\in\H$, $y\in Y$ and $z\in Z$.
Theorem 3.2 in \cite{EvSo1984} implies easily the following result:

\begin{remark}\label{euclidean assumption}
Let us consider the problem (\ref{zero sum}) with the assumptions 1., 2'.  and 3'..
 Then $V^-$ is bounded and uniformly
Lipschitz continuous w.r.t. the Euclidean distance $d_E$, i.e.
$$|V^-(t,x)-V^-(t',x')|\le C(|t-t'|+\|x-x'\|),$$
for every $t,t'\in [0,T]$ and $x,x'\in\H$. Consequently, $V^-$ is
Lipschitz continuous w.r.t. the gauge distance $d_G$. A similar result holds for $V^+$.
\end{remark}

\noindent We note that our result in Theorem \ref{Teo visco} is more precise under weaker assumptions, since there exists $d_G$-Lipschitz functions that are not $d_E$-Lipschitz.

\medskip
 \noindent
\textbf{Proof of Theorem \ref{Teo visco}}  The idea of the proof follows from Theorem 3.2 in \cite{EvSo1984}, but here we use all the fine properties of the horizontal curves in $\H$ w.r.t. the $d_G$-distance proved in subsection \ref{curve orizzonatali}.

Let us fix $\tau<\tau'$ in $[0,T]$ and $\xi,\xi'\in\H$.
It is immediate to see that
$$
|V^-(\tau,\xi)|\le C_1T+C_2.
$$
Now let us fix $\epsilon>0$. There exists $\widehat\beta\in {\mathcal S}_\beta(\tau)$ such that
\begin{equation}\label{dim2a}
V^-(\tau,\xi)\ge \sup_{y\in
{\mathcal Y}(\tau)}\left\{ \int_\tau^T F(t,x,y,\widehat \beta[y])\de
t+g(x(T)) \right\} -\epsilon.
\end{equation}
Fix $y_0\in Y$. For every $y\in
{\mathcal Y}(\tau')$, let us define $\widehat y\in
{\mathcal Y}(\tau) $ by
\begin{equation}\label{dim2aa}
\widehat y(t)=\left\{
              \begin{array}{ll}
                 y_0, & \hbox{for $t\in [\tau,\tau')$} \\
                 y(t), & \hbox{for $t\in [\tau',T]$}
              \end{array}
            \right.
\end{equation}
Let us define $\widetilde\beta\in {\mathcal S}_\beta(\tau')$ such that
$$
\widetilde\beta[y](t)=\widehat\beta[\widehat y](t)
,\qquad \forall y\in  {\mathcal Y}(\tau'),\ t\in [\tau',T]$$
Clearly,
$$
V^-(\tau',\xi')\le \sup_{y\in
{\mathcal Y}(\tau')}\left\{ \int_{\tau'}^T F(t,x,y,\widetilde \beta[y])\de
t+g(x(T)) \right\}.
$$
Let $\widetilde y\in {\mathcal Y}(\tau')$ be such that
\begin{equation}\label{dim2b}
V^-(\tau',\xi')\le  \int_{\tau'}^T F(t,x,\widetilde y,\widetilde \beta[\widetilde y])\de
t+g(x(T))+\epsilon.
\end{equation}
From \eqref{dim2a} we get
\begin{equation}\label{dim2c}
V^-(\tau,\xi)\ge \int_\tau^T F(t,x,\widehat y,\widehat \beta[\widehat y])\de
t+g(x(T)) -\epsilon,
\end{equation}
where $\widehat y$ is defined by $\widetilde y$ via  \eqref{dim2aa} replacing $y$ with  $\widetilde y$.
Note that the trajectories $x$ that appear in \eqref{dim2b} and in \eqref{dim2c} are different functions. In particular,
denoting by $\widetilde  x$ and $\widehat x$  such trajectories in  \eqref{dim2b} and in \eqref{dim2c} respectively, we have that
$\widetilde  x$ is a horizontal curve on $[\tau',T]$ with horizontal velocity $-\widetilde \beta[\widetilde y]$ and initial point $\xi'$,  while
$\widehat  x$ is a horizontal curve on $[\tau,T]$ with horizontal velocity $-\widehat \beta[\widehat y]$ and initial point $\xi$.
Since  $\widetilde y=\widehat y$ and $\widetilde \beta[\widetilde y]= \widehat \beta[\widehat y]$ on $[\tau',T]$,  by Proposition \ref{hor curve new 2} it is easy to prove that, for every $t\in[\tau',T]$,
$$
d_G\left(\widetilde x(t),\widehat x(t)\right)
 \le \widetilde C \left(d_G(\xi',\xi)+(\tau'-\tau)\right).
$$
By \eqref{dim2b} and \eqref{dim2c}, assumptions \emph{2.} and  \emph{3.}
we obtain
\begin{eqnarray}
V^-(\tau',\xi')
-
V^-(\tau,\xi)
&\le&  \int_{\tau'}^T \left( F(t,\widetilde x,\widetilde y,\widetilde \beta[\widetilde y])- F(t,\widehat x,\widehat y,\widehat \beta[\widehat y])\right)\de t+\nonumber\\
&&\quad\quad -\int_\tau^{\tau'} F(t,\widehat x,\widehat y,\widehat \beta[\widehat y])\de
t+g(\widetilde x(T))-g(\widehat x(T)) +2\epsilon\nonumber\\
&\le&  C_1' \int_{\tau'}^T d_G(\widetilde x(t),\widehat x(t))\de
t+
(\tau'-\tau)C_1+ C_2'd_G\left(\widetilde x(T),\widehat x(T)\right)
+2\epsilon\nonumber\\
&\le&   C'\left( d_G(\xi',\xi)+(\tau'-\tau)\right)
+2\epsilon\label{dim2d}
\end{eqnarray}
with
\begin{equation}\label{C prime}
 C'=\widetilde C(C_1'T+C_2')+C_1.
\end{equation}
This concludes the first part of the proof.

Let $\epsilon$ again be fixed. Then there exists $\widehat\beta\in {\mathcal S}_\beta(\tau')$ such that
\begin{equation}\label{dim2abis}
V^-(\tau',\xi')\ge \sup_{y\in
{\mathcal Y}(\tau')}\left\{ \int_{\tau'}^T F(t,x,y,\widehat \beta[y])\de
t+g(x(T)) \right\} -\epsilon.
\end{equation}
For every $y\in
{\mathcal Y}(\tau)$, let us define $\widehat y\in
{\mathcal Y}(\tau') $ by
\begin{equation}\label{dim2aabis}
\widehat y(t)=y(t),\qquad  \forall t\in [\tau',T]
\end{equation}
Fix $y_0\in Y$. Let us define $\widetilde\beta\in {\mathcal S}_\beta(\tau)$ such that, for every $y\in
{\mathcal Y}(\tau)$
$$
\widetilde \beta[y](t)=\left\{
              \begin{array}{ll}
                 y_0, & \hbox{for $t\in [\tau,\tau')$} \\
                \widehat \beta[\widehat y](t), & \hbox{for $t\in [\tau',T]$}
              \end{array}
            \right.
$$
Clearly,
$$
V^-(\tau,\xi)\le \sup_{y\in
{\mathcal Y}(\tau)}\left\{ \int_{\tau}^T F(t,x,y,\widetilde \beta[y])\de
t+g(x(T)) \right\}.
$$
Let $\widetilde y\in {\mathcal Y}(\tau)$ be such that
\begin{equation}\label{dim2bbis}
V^-(\tau,\xi)\le  \int_{\tau}^T F(t,\widetilde x,\widetilde y,\widetilde \beta[\widetilde y])\de
t+g(x(T))+\epsilon.
\end{equation}
The inequality \eqref{dim2abis} gives
\begin{equation}\label{dim2cbis}
V^-(\tau',\xi')\ge \int_{\tau'}^T F(t,\widehat x,\widehat y,\widehat \beta[\widehat y])\de
t+g(x(T)) -\epsilon.
\end{equation}
where $\widehat y$ is defined by $\widetilde y$ via relation \eqref{dim2aabis} replacing $y$ with $\widetilde y$. Note that
 $\widetilde  x$  is a horizontal curve on $[\tau,T]$ with horizontal velocity $-\widetilde \beta[\widetilde y]$ and initial point $\xi$,  while
$\widehat  x$ is a horizontal curve on $[\tau',T]$ with horizontal velocity $-\widehat \beta[\widehat y]$ and initial point $\xi'$.
Since  $\widetilde y=\widehat y$ and $\widetilde \beta[\widetilde y]= \widehat \beta[\widehat y]$ on $[\tau',T]$,
 by
 Proposition  \ref{hor curve new 2},
we have that, for every $t\in[\tau',T]$,
 $$
d_G\left(\widehat x(t),\widetilde x(t)\right)
 \le
\widetilde C \left(d_G(\xi',\xi)+(\tau'-\tau)\right).
$$
By \eqref{dim2bbis} and \eqref{dim2cbis}, assumptions \emph{2.} and  \emph{3.}
we obtain
\begin{eqnarray}
V^-(\tau,\xi)
-
V^-(\tau',\xi')
&\le&  \int_{\tau'}^T \left( F(t,\widetilde x,\widetilde y,\widetilde \beta[\widetilde y])- F(t,\widehat x,\widehat y,\widehat \beta[\widehat y])\right)\de
t+\nonumber\\
&&\quad\quad +\int_\tau^{\tau'} F(t,\widetilde x,\widetilde y,\widetilde \beta[\widetilde y])\de
t+g(\widetilde x(T))-g(\widehat x(T))
+2\epsilon\nonumber\\
&\le&  C_1' \int_{\tau'}^T d_G(\widetilde x(t),\widehat x(t))\de
t+
(\tau'-\tau)C_1+ C_2'd_G(\widetilde x(T),\widehat x(T))
+2\epsilon\nonumber\\
&\le&  C'\left( d_G(\xi',\xi)+(\tau'-\tau)\right)
+2\epsilon\label{dim2d bis}
\end{eqnarray}
with $ C'$ as in  \eqref{C prime}.
This inequality and \eqref{dim2d} conclude the proof.
\finedim

In the fundamental paper \cite{Pa1989}, Pansu provides a Rademacher--Stefanov type result in the Carnot group setting; in particular, he proves that every Lipschitz continuous function w.r.t. a homogeneous distance on $\H$ is differentiable almost everywhere in the horizontal directions.
Hence, our previous result implies that the lower value function admits the horizontal gradient and the derivative w.r.t.  $t$, i.e.
$$\nabla_HV^-(t,x)=\left(X_1V^-(t,x),X_2V^-(t,x)\right)\qquad \texttt{\rm and}\quad \frac{\partial V^-}{\partial t}(t,x),$$ for almost everywhere $(t,x)\in[0,T]\times \H$.
Therefore, $V^-$ could be a candidate for a viscosity solution, as we will see in definition \ref{def vis sol}.

\noindent
Moreover, a more precise estimate in  \eqref{dim2d} and in \eqref{dim2d bis} gives
$$|
V^-(t,x)
-
V^-(t',x')|\le  \widetilde C(C_1'T+C_2')\left( d_G(x',x)+|t'-t|\right)+C_1|t'-t|
$$
This implies that, taking into account (\ref{C tilde 0}), we have the following $d_G$-Lipschitz constant
\begin{remark}\label{rem C sharp}
For a.e. $(t,x)\in[0,T]\times \H$ we have
$$
\|\nabla_H V^-(t,x)\|\le (1+3R_Z)e^{\frac{TR_Z }{2}}(C_1'T+C_2'):=C^\sharp
$$
\end{remark}
It is important to notice that $C^\sharp$ does not depend on $R_Y$.

\section{Viscosity solutions for Hamilton--Jacobi--Isaacs equation}

This section is essentially devoted to the proof of Theorem \ref{Teo vis2}.
In order to recall the notion of viscosity solution in our context (see for example \cite{MaSt2002}) we
  say that, given an open interval $I$, a function $\psi:I\times\H\to\R$  is in $\Gamma^1(I\times\H)$ if $(t,x)\mapsto\left(\frac{\partial \psi(t,x)}{\partial t},\ X_1\psi(t,x),\ X_2\psi(t,x)\right)$ is a continuous function.

\begin{definition}[viscosity solution]\label{def vis sol} Let $\Hcal:[0,T]\times\H\times\R^{2}\to \R$ be a continuous function and let $u:[0,T]\times\H \to \R$ be a bounded and uniformly continuous function, with $u(T,x)=g(x)$ in $\H$. We say that $u$  is a viscosity subsolution of the Hamilton--Jacobi equation
\begin{equation}\label{Ham Jac eq}
\left\{
  \begin{array}{ll}
    \displaystyle \frac{\partial u}{\partial t}(t,x)+\Hcal(t,x,\nabla_{H}u(t,x))=0 & \hbox{in $(0,T)\times\H$} \\
    u(T,x)=g(x) & \hbox{in $\H$}
  \end{array}
\right.
\end{equation}
if, whenever $(t_0,x_0)\in (0,T)\times\H$ and $\psi$ is a test function in $\Gamma^1((0,T)\times\H)$
touching $u$ from above at $(t_0,x_0)$, i.e.
$$
u(t_0,x_0)=\psi(t_0,x_0)\quad {\sl and}\ u(t,x)\le\psi(t,x) \quad \texttt{\sl in a neighborhood of}\ (t_0,x_0),
$$
we have
\begin{equation}\label{dis vis sol 1}
\frac{\partial \psi}{\partial t}(t_0,x_0)+\Hcal(t_0,x_0,\nabla_{H}\psi(t_0,x_0))\ge 0.
\end{equation}
We say that $u$ is a viscosity supersolution of equation \eqref{Ham Jac eq}
if, whenever $(t_0,x_0)\in (0,T)\times\H$ and $\psi$ is a  test function in $\Gamma^1((0,T)\times\H)$
touching $u$ from below at $(t_0,x_0)$, i.e.
$$
u(t_0,x_0)=\psi(t_0,x_0)\quad {\it and}\ u(t,x)\ge\psi(t,x) \quad \texttt{\it in a neighborhood of}\ (t_0,x_0),
$$
we have
\begin{equation}\label{dis vis sol 2}
\frac{\partial \psi}{\partial t}(t_0,x_0)+\Hcal(t_0,x_0,\nabla_{H}\psi(t_0,x_0))\le 0.
\end{equation}
A  function that is both a viscosity subsolution and a viscosity supersolution is called viscosity solution.
\end{definition}
An equivalent definition  of viscosity solution \eqref{Ham Jac eq} uses the notion of supejets (see \cite{MaSt2002}, \cite{Ma2003}).
We recall that if in \eqref{Ham Jac eq} we change the condition $ u(T,x)=\psi(x)$ with an initial condition of the type
$$
 u(0,x)=\psi(x)\qquad x\in\H,
$$
then the viscosity solution of the new problem is defined by reversing the inequalities \eqref{dis vis sol 1} and \eqref{dis vis sol 2}.

\begin{definition}[upper and lower Hamiltonian]\label{def upper lower ham}
Let us consider the zero game \eqref{zero sum}.
We define the {lower Hamiltonian}
 $H^-:[0,T]\times\H\times\R^{2}\to\R$  by
\begin{equation}\label{def H-}
H^-(t,x,\lambda)=\max_{y\in Y}\min_{z\in Z}\Bigl( F(t,x,y,z)-\lambda\cdot z\Bigr)
\end{equation}
and
 the {upper Hamiltonian }
 $H^+:[0,T]\times\H\times\R^{2}\to\R$  by
$$
H^+(t,x,\lambda)=\min_{z\in Z}\max_{y\in Y}\Bigl( F(t,x,y,z)-\lambda\cdot z\Bigr).
$$
\end{definition}
It is easy to see that $H^-\le H^+$.
We say that the $\min\max$ condition, or Isaacs' condition, is satisfied if
$H^− = H^+.$
 In this case, we define the {Hamiltonian} $H$ by
$$
H (t, x,\lambda) = H^-(t, x,\lambda) = H^+(t, x,\lambda).$$

Let us spend few lines to make some comments on the Definition \ref{def upper lower ham}.
In a classical zero game case, if we have a trajectory $x$ in $\R^n$,
then we usually introduce a multiplier $\lambda$ with the same dimension, i.e. $\lambda\in\R^n$; more precisely, if $\dot{x}=h(t,x,y,z)$ is the dynamics of the zero game, in the definition \eqref{def H-} of $H^-$  the function $(F+\l\cdot g)$ appears as argument of the $\max\min$. In our case,
 the trajectory $x$ is in $\H$ but the multiplier $\l$ is 2-dimensional and  takes into account only the horizontal velocity of the horizontal curve $x$.

Now we are ready to prove Theorem \ref{Teo vis2 intro}, i.e.

\begin{theorem}\label{Teo vis2}
Let us consider the problem (\ref{zero sum}) with the assumptions
1., 2. and 3.. Then $V^-$ is a viscosity solution
of the lower Hamilton--Jacobi--Isaacs equation
 \begin{equation}\label{system PD-}
 \left\{
   \begin{array}{ll}
   {\displaystyle \frac{\partial u}{\partial t}(t,x)+H^-(t,x,\nabla_{H}u(t,x))=0} & \hbox{for $(t,x)\in(0,T)\times\H$} \\
  u(T,x)=g(x)&  \hbox{for
$x\in\H$}
   \end{array}
 \right.
 \end{equation}
 and $V^+$ is a  viscosity solution
of the upper Hamilton--Jacobi--Isaacs equation
 \begin{equation}\label{system PD+}
 \left\{
   \begin{array}{ll}
   {\displaystyle \frac{\partial u}{\partial t}(t,x)+H^+(t,x,\nabla_{H}u(t,x))=0} & \hbox{for $(t,x)\in(0,T)\times\H$} \\
  u(T,x)=g(x)&  \hbox{for
$x\in\H$}
   \end{array}
 \right.
 \end{equation}
  \end{theorem}

\noindent
The proof of this theorem requires the following lemma:
\begin{lemma}\label{lemma}
Let $\psi\in\Gamma^1((0,T)\times H)$.


If there exists $\theta>0$ such that
\begin{equation}\label{dim vis sol 2}
\frac{\partial \psi}{\partial t}(t_0,x_0)+H^-(t_0,x_0,\nabla_{H}\psi(t_0,x_0))\ge\theta,
\end{equation}
then,  for all sufficiently small $\tau>0$,
 there exists   $\widetilde y\in  {\mathcal Y}(t_0)$  such that for every $\widetilde\beta
\in {\mathcal S}_\beta(t_0)$ we have
\begin{equation}\label{lemma a}
\int_{t_0}^{t_0+\tau}\left(
\frac{\partial \psi}{\partial t}(s,\widetilde x(s))+F(s,\widetilde x(s),\widetilde y(s),\widetilde\beta[\widetilde y](s))- \widetilde\beta[\widetilde y](s)\cdot\nabla_H\psi(s,\widetilde x(s))\right)d s\ge\frac{\tau\theta}{2},
\end{equation}
where $\widetilde x$ is the horizontal curve on $[t_0,t_0+\tau]$ with horizontal velocity $\widetilde\beta[\widetilde y]$ and initial point $x_0$.

If there exists $\theta>0$ such that
\begin{equation}\label{dim vis sol 2 bis}
\frac{\partial \psi}{\partial t}(t_0,x_0)+H^-(t_0,x_0,\nabla_{H}\psi(t_0,x_0))\le-\theta,
\end{equation}
then,  for all sufficiently small $\tau>0$,
 there exists   $\widetilde\beta
\in {\mathcal S}_\beta(t_0)$ such that for every $\widetilde y\in  {\mathcal Y}(t_0)$   we have
\begin{equation}\label{lemma b}
\int_{t_0}^{t_0+\tau}\left(
\frac{\partial \psi}{\partial t}(s,\widetilde x(s))+F(s,\widetilde x(s),\widetilde y(s),\widetilde\beta[\widetilde y](s))-\widetilde\beta[\widetilde y](s)\cdot\nabla_H\psi(s,\widetilde x(s))\right)d s\le-\frac{\tau\theta}{2},
\end{equation}
where $\widetilde x$ is as before.


\end{lemma}
\noindent
The proof of this lemma can be found in \cite{EvSo1984} (see Lemma 4.3, where a classical gradient instead of our horizontal gradient appears) and it is based on the continuity of the function
$$(t,x)\mapsto \frac{\partial \psi}{\partial t}(t,x)+F(t,x,y,z)-z\cdot \nabla_{H}\psi(t,x)$$
and on the compactness of the control sets $Y$ and $Z$.

\noindent
\textbf{Proof of Theorem \ref{Teo vis2}} The proof follows the idea of Theorem 4.1 in \cite{EvSo1984} and uses the properties of the horizontal curves in subsection  \ref{curve orizzonatali}.
It is obvious, by definition, that $V^-(T,x)=g(x)$, for every
$x\in\H$. So, let us fix $(t_0,x_0)\in(0,T)\times\H$.

\noindent First, let $\psi\in \Gamma^1((0,T)\times\H)$  be a
test function
touching $V^-$ from below at $(t_0,x_0)$, i.e.
\begin{equation}\label{dim vis sol 1}
V^-(t_0,x_0)=\psi(t_0,x_0)\quad {\rm and}\ V^-(t,x)\ge\psi(t,x) \quad \texttt{\rm in a neighborhood of}\ (t_0,x_0).
\end{equation}
We have to prove that (\ref{dis vis sol 2}) holds with $\Hcal=H^-.$
By contradiction, let us assume that this is not true and that there exists $\theta>0$ such that holds (\ref{dim vis sol 2}); then, \eqref{lemma a} implies
that
\begin{equation}\label{dim vis sol 4a}
\inf_{\beta\in {\mathcal S}_\beta(t_0)}\sup_{y\in
{\mathcal Y}(t_0)}
\int_{t_0}^{t_0+\tau}\left(
\frac{\partial \psi}{\partial t}(s, x)+F(s, x, y,\beta[y])- \beta[ y]\cdot\nabla_H\psi(s, x)\right)d s\ge\frac{\tau\theta}{2}
\end{equation}
where
 $ x$ solves
\begin{equation}\label{dim vis sol 5}
\left\{
    \begin{array}{ll}
      \dot{ x}=-f^\H  (x,\beta[ y])& \hbox{in $[t_0,T]$} \\
       x(t_0)=x_0, &     \end{array}
  \right.
\end{equation}
Now by Theorem \ref{Teo vis0} we know that
\begin{equation}\label{dim vis sol 6}
V^-(t_0,x_0)=\inf_{\beta\in {\mathcal S}_\beta(t_0)}\sup_{y\in
{\mathcal Y}(t_0)}\left\{ \int_{t_0}^{t_0+\tau} F(s,x,y,\beta[y])\de s+V^-(t_0+\tau,x(t_0+\tau)) \right\}
\end{equation}
with $x$ as before. For every such horizontal curve $x$, \eqref{dim vis sol 1} and Proposition \ref{bounded traj} imply that, for $\tau$ small enough,
\begin{equation}\label{dim vis sol 7}
0=V^-(t_0,x_0)-\psi(t_0,x_0)\le V^-(t_0+\tau,x(t_0+\tau))-\psi(t_0+\tau,x(t_0+\tau)) \end{equation}
Since $x$ is horizontal and $\psi$ is in $\Gamma^1$, Remark \ref{derivata} and \eqref{dim vis sol 5} imply
\begin{eqnarray}\psi(t_0+\tau,x(t_0+\tau)) -\psi(t_0,x_0)&=&\int_{t_0}^{t_0+\tau}
\frac{d\psi (s,x(s))}{d s} ds\nonumber\\
&=&\int_{t_0}^{t_0+\tau} \left(\frac{\partial \psi}{\partial t}(s, x(s))-\beta[y](s)\cdot \nabla_H \psi(s,x(s))\right) ds\qquad\label{dim vis sol 8}
\end{eqnarray}
Relations \eqref{dim vis sol 6}--\eqref{dim vis sol 8}  give
\begin{eqnarray*}
0&\ge& \inf_{\beta\in {\mathcal S}_\beta(t_0)}\sup_{y\in
{\mathcal Y}(t_0)}\left\{ \int_{t_0}^{t_0+\tau} F(s,x,y,\beta[y])\de s+\psi(t_0+\tau,x(t_0+\tau))-\psi(t_0,x_0) \right\}\\
&=& \inf_{\beta\in {\mathcal S}_\beta(t_0)}\sup_{y\in
{\mathcal Y}(t_0)}\left\{ \int_{t_0}^{t_0+\tau}
\left(F(s,x,y,\beta[y])+\frac{\partial \psi}{\partial t}(s, x)-\beta[y]\cdot \nabla_H \psi(s,x)\right)\de s\right\}\end{eqnarray*}
This inequality contradicts \eqref{dim vis sol 4a}, hence  \eqref{dim vis sol 2} is false and this concludes the first part of the proof.

Now, let $\psi\in \Gamma^1((0,T)\times\H)$   be a
test function
touching $V^-$ from above at $(t_0,x_0)$, i.e.
\begin{equation}\label{dim vis sol 1 bis}
V^-(t_0,x_0)=\psi(t_0,x_0)\quad {\rm and}\ V^-(t,x)\le\psi(t,x) \quad \texttt{\rm in a neighborhood of}\ (t_0,x_0).
\end{equation}
We have to prove that (\ref{dis vis sol 1}) holds with $\Hcal=H^-.$
Let us assume that this is not true and that there exists $\theta>0$ such that
(\ref{dim vis sol 2 bis}) holds; then \eqref{lemma b}  implies
that
\begin{equation}\label{dim vis sol 4a bis}
\inf_{\beta\in {\mathcal S}_\beta(t_0)}\sup_{y\in
{\mathcal Y}(t_0)}
\int_{t_0}^{t_0+\tau}\left(
\frac{\partial \psi}{\partial t}(s, x)+F(s, x, y,\beta[y])- \beta[ y]\cdot\nabla_H\psi(s, x)\right)d s\le-\frac{\tau\theta}{2}
\end{equation}
where
 $ x$ is as in \eqref{dim vis sol 5}.
 For every such horizonal curve $x$, requirement \eqref{dim vis sol 1 bis} and Proposition \ref{bounded traj} imply that, for $\tau$ small enough,
\begin{equation}\label{dim vis sol 7 bis}
0=V^-(t_0,x_0)-\psi(t_0,x_0)\ge V^-(t_0+\tau,x(t_0+\tau))-\psi(t_0+\tau,x(t_0+\tau)) \end{equation}
Relations \eqref{dim vis sol 6}, \eqref{dim vis sol 8} and \eqref{dim vis sol 7 bis}  give
\begin{eqnarray*}
0&\le& \inf_{\beta\in {\mathcal S}_\beta(t_0)}\sup_{y\in
{\mathcal Y}(t_0)}\left\{ \int_{t_0}^{t_0+\tau} F(s,x,y,\beta[y])\de s+\psi(t_0+\tau,x(t_0+\tau))-\psi(t_0,x_0) \right\}\\
&=& \inf_{\beta\in {\mathcal S}_\beta(t_0)}\sup_{y\in
{\mathcal Y}(t_0)}\left\{ \int_{t_0}^{t_0+\tau}
\left(F(s,x,y,\beta[y])+\frac{\partial \psi}{\partial t}(s, x)-\beta[y]\cdot \nabla_H \psi(s,x)\right)\de s\right\}.\end{eqnarray*}
This inequality contradicts \eqref{dim vis sol 4a bis}: hence  \eqref{dim vis sol 2 bis} is false and this concludes the proof for $V^-$.
 In a similar way one proves that $V^+$ is a viscosity solution of \eqref{system PD+}.
\finedim

\section{Representation of solutions of Hamilton--Jacobi equations}

We are now in the position to study the viscosity solution of our initial Hamilton--Jacobi problem  \eqref{Ham Jac eq intro}, i.e.
 \begin{equation}\label{pro finale}
 \left\{
   \begin{array}{ll}
   {\displaystyle \frac{\partial u}{\partial t}(t,x)+\Hcal(t,x,\nabla_H u(t,x))=0} & \hbox{for $(t,x)\in(0,T)\times\H$} \\
  u(0,x)=g(x)&  \hbox{for
$x\in\H$}
   \end{array}
 \right.
 \end{equation}
under assumptions \emph{3.} and \emph{4.}.

Having in mind problem \eqref{pro finale}, let us introduce a zero sum game as follows:  set
\begin{equation}\label{costanti finale}
R_Z=K,\qquad R_Y=(1+3K)e^{\frac{TK }{2}}(D_1'T+C_2')
\end{equation}
in assumption \emph{1.}, and consider the function
 \begin{equation}\label{F nuova}
F(t,x,y,z)=-\Hcal(T-t,x,y)+z\cdot y. \end{equation}
Relations \eqref{costanti finale} and \eqref{F nuova} give us a zero sum game as in \eqref{zero sum} associated to our initial problem \eqref{pro finale}.
It is clear that  assumptions \emph{1.}  and \emph{4.}  guarantee that $F$ in \eqref{F nuova} satisfies assumption \emph{2.} with $C_1=D_1+R_ZR_Y$ e $C_1'=D_1'.$ Hence, Theorem \ref{Teo visco} implies that the lower value function $V^-_F$ for the zero game  \eqref{zero sum}, with $F$ as in \eqref{F nuova} and $R_Y$ and $R_Z$ as in \eqref{costanti finale},
 \begin{eqnarray*}
V^-_F(\tau,\xi)&=&\inf_{\beta\in {\mathcal S}_\beta(\tau)}\sup_{y\in
{\mathcal Y}(\tau)}\left\{ \int_\tau^T \Bigl(-\Hcal(T-t,x,y)+\beta[y]\cdot y\Bigr)\de
t+g(x(T)) \right\}\\
&&\texttt{\rm with}\ x(t)=\xi-\int_\tau^t f^\H(x,\beta[y])\de s,
\end{eqnarray*}
is bounded and $d_G$-Lipschitz w.r.t. $x$. Remark \ref{rem C sharp} gives
$\|\nabla_H V_F^-(t,x)\|\le C^\sharp$, with
 $$
 C^\sharp=(1+3K)e^{\frac{TK }{2}}(D_1'T+C_2').
 $$
Note that $C^\sharp=R_Y$. Moreover, by Theorem \ref{Teo vis2}, $V^-_F$ is a viscosity solution of the lower Hamilton--Jacobi--Isaacs equation \eqref{system PD-}, i.e.
\begin{equation}\label{equa Isa Ham Jac fin}
 \left\{
   \begin{array}{ll}
   {\displaystyle \frac{\partial u}{\partial t}(t,x)+H^-(t,x,\nabla_{H}u(t,x))=0} & \hbox{for $(t,x)\in(0,T)\times\H$} \\
  u(T,x)=g(x)&  \hbox{for
$x\in\H$}
   \end{array}
 \right.
\end{equation}
where, as in \eqref{def H-},
\begin{equation}\label{H finale}
\displaystyle H^-(t,x,\l)=\max_{y\in Y}\min_{z\in Z}\left(- \Hcal(T-t,x,y)+z\cdot y-\l\cdot z\right),
\end{equation}
and $Y=\overline{B_{\R^2}(0,R_Y)}$, $Z=\overline{B_{\R^2}(0,R_Z)}$.
Clearly, assumption \emph{4.} implies
$$
\Hcal(t,x,\l)\le \Hcal(t,x,y)+K\|\l-y\|
$$
for every $y,\l\in Y$. Hence, taking into account that $K=R_Z$, for every $t\in [0,T],\ x\in\H$ and $\l\in Y$,
 \begin{eqnarray*}
\Hcal(t,x,\l)&=&\min_{y\in Y}\left( \Hcal(t,x,y)+R_Z\|\l-y\|\right)\\
&=&\min_{y\in Y}\max_{z\in Z}\left( \Hcal(t,x,y)+z\cdot(\l-y)\right)\\
&=&-\max_{y\in Y}\min_{z\in Z}\left(- \Hcal(t,x,y)-z\cdot(\l-y)\right).
\end{eqnarray*}
The equality above and \eqref{H finale} imply that
\begin{equation}\label{H uguali}
H^-(T-t,x,\l)=-\Hcal(t,x,\l)\qquad \forall t\in[0,T],\ x\in\H,\ \l\in Y.
\end{equation}
The function $U$, defined by $U(t,x)=V^-_F(T-t,x)$ and \eqref{equa Isa Ham Jac fin},
is a viscosity solution for
$$
 \left\{
   \begin{array}{ll}
   {\displaystyle \frac{\partial u}{\partial t}(t,x)-H^-\left(T-t, x,\nabla_H u(t,x)\right)=0} & \hbox{for $(t,x)\in(0,T)\times\H$} \\
  u(0,x)=g(x)&  \hbox{for
$x\in\H$}
   \end{array}
 \right.
 $$
Taking into account that $\|\nabla_H U(t,x)\|\le C^\sharp=R_Y$, and using \eqref{H uguali}, we finally have the following representation for the viscosity solution for the Hamilton--Jacobi equation \eqref{pro finale}.

\begin{theorem}\label{ciao}
Let us consider  problem \eqref{pro finale} with the assumptions  1. ($R_Z$ and $R_Y$ as in \eqref{costanti finale}), 3. and 4.
Then, the function $U$ defined by
 \begin{eqnarray*}
&&U(\tau,\xi)= \inf_{\beta\in {\mathcal S}_\beta(T-\tau)}\sup_{y\in
{\mathcal Y}(T-\tau)}\left\{ \int_{T-\tau}^T \Bigl(-\Hcal(T-t,x,y)+\beta[y]\cdot y\Bigr)\de
t+g(x(T)) \right\}
\\
&&\qquad\qquad\texttt{\rm with}\ x(t)=\xi-\int_{T-\tau}^t f^\H(x,\beta[y])\de s,\nonumber
\end{eqnarray*}
is $d_G$-Lipschitz and is a viscosity solution  for \eqref{pro finale}.

  \end{theorem}

\subsection{A particular case and a question.}
Let us consider the particular case $\Hca=\Hca(y)$. The previous arguments give that, under  the same assumptions  \emph{1.}, \emph{3.} and \emph{4.},
 the function $U$ defined by
 \begin{eqnarray}
&&U(\tau,\xi)= \inf_{\beta\in {\mathcal S}_\beta(0)}\sup_{y\in
{\mathcal Y}(0)}\left\{ \int_0^\tau\Bigl(-\Hcal(y)+\beta[y]\cdot y\Bigr)\de
t+g(x(\tau)) \right\}\label{fun U}
\\
&&\qquad\qquad\texttt{\rm with}\ x(t)=\xi-\int_{0}^\tau f^\H(x,\beta[y])\de s,\nonumber
\end{eqnarray}
with $R_Z$ and $R_Y$ as in \eqref{costanti finale}, is a viscosity solution  for
\begin{equation}\label{pro finale autonomo}
 \left\{
   \begin{array}{ll}
   {\displaystyle \frac{\partial u}{\partial t}(t,x)+\Hcal\left(\nabla_H u(t,x)\right)=0} & \hbox{for $(t,x)\in(0,T)\times\H$} \\
  u(0,x)=g(x)&  \hbox{for
$x\in\H$.}
   \end{array}
 \right.
 \end{equation}
 Theorem 4 in \cite{MaSt2002} guarantees the uniqueness of such solution:
\begin{proposition}
In the assumptions 1.,  3. and 4., the function $U$ in \eqref{fun U} is $d_G$-Lipschitz and is the unique viscosity solution of \eqref{pro finale autonomo} satisfying
$$
\lim_{t\to 0}\sup_{x\in \H}\left|u(x,t)-u(x,0)\right|=0.
$$
\end{proposition}
It is well known that, if the Hamiltonian function depends only on the gradient, it is possible to write a Hopf--Lax formula: the first result in this line of investigation in $\H$ is in \cite{MaSt2002}, a more general result can be found in \cite{BaCaPi2014}. Hence, a very interesting question is the following: if we add to the previous assumptions the following
\begin{itemize}
\item[\emph{5.}] the function  $g:\H\to\R$ is \lq\lq convex\rq\rq,
\end{itemize}
is it possible to apply the ideas in \cite{BaEv1984} to obtain a Hopf--Lax formula for the function $U$ in  \eqref{fun U}?
If $g$ is convex in the $\R^3$ classical sense, and hence is locally $d_E$-Lipschitz, one can try to apply the same arguments of section 3 in \cite{BaEv1984} where a Jensen inequality plays a fundamental role.

But the very interesting and natural question, taking into account the Sub--Riemannian setting and  assumption \emph{4.}, arises if we require that $g$ is only H--convex, i.e.,  for every fixed  $x\in\H$ and $w\in V_1$ the function   $s\mapsto g(x\circ\exp(sw))$ is convex. The simplest example of a H--convex function, but not $\R^3$--convex, is  $x\mapsto\|x\|_G$; moreover, we recall that an H--convex function is $d_G$-Lipschitz  (see \cite{CaPi2011} for details on the properties of these H--convex functions). Unfortunately with this notion of convexity, to our knowledge, there is not in the literature a Jensen--type inequality in the Heisenberg group that would allow to follow the ideas in  \cite{BaEv1984} in order to obtain a Hopf--Lax formula: this is a very big obstacle. We are working on this obstacle.

\end{document}